# An algorithm for minimum cardinality generators of cones

**Matthias G. Mayer***
matthias.georg.mayer@gmail.com

**Fabian von der Warth***
vdwarth@rptu.de

**ABSTRACT**

This paper presents a novel proof that for any convex cone, the size of conically independent generators is at most twice that of minimum cardinality generators. While this result is known for linear spaces, we extend it to general cones through a decomposition into linear and pointed components. Our constructive approach leads to a polynomial-time algorithm for computing minimum cardinality generators of finitely generated cones, improving upon existing methods that only compute conically independent generators.

**Keywords:** conically indepedent; cone; frame; minimum cardinality
**MSC Subject Classification:** 90C25

## 1 Introduction

A cone is a subset of $\mathbb{R}^n$ that is closed under positive linear combinations. A generator of a cone $C$ is a subset of $\mathbb{R}^n$ whose conic hull, i.e. the smallest cone that contains $S$, is $C$. A conically independent generator is a generator for which no proper subset is still a generator of the same cone. For formal definitions we refer to Section 2.

For cones that are linear spaces, it is well known that conically independent generators may exceed minimum cardinality generators in size by a factor of at most two [1, Theorem 6.7]. Although these fundamental theorems are well-established, they are typically presented as parts of broader mathematical treatises [1]. To the best of our knowledge, no comprehensive treatment of general cones, as presented here, exists in the literature. For finitely generated cones, finding minimum cardinality generator representations remains less well understood. [2] provides only an algorithm for computing conically independent generators. We improve on [2] by providing a polynomial-time algorithm for finding minimum cardinality generators of cones. We establish this result within a broader theorem that demonstrates that for any cone, the size of conically independent generators is at most twice that of minimum cardinality generators. Our proof is self-contained, except for an elementary result of Dattoro [3] (cf. Theorem 8), and provides an alternative approach to that of [1].

This finds application for instance in the context of objective reduction in multi-objective optimization. In [4], the notion of redundant objective functions is defined within linear problems. An objective function is said to be redundant if its removal does not change the set of efficient solutions. The condition that $F_i = \sum_{j \neq i} \lambda_j F_j$ with $\lambda_j \geq 0, j \neq i$ is proven to be sufficient

---

for an objective function to be redundant. This is equivalent to the fact that $\{F_i : i \in I\}$ is not conically independent. Since any minimum cardinality generator is also conically independent, it is advantageous to determine a minimum cardinality generator within this context.

## 2 Preliminaries

For the rest of the section we let $S \subseteq \mathbb{R}^n$ be finite, for some $n \in \mathbb{N}_{\geq 1}$.

**Definition 1** (Convex cone, generation): A subset $C \subseteq \mathbb{R}^n$ is called a convex cone (or just cone) if it is closed under positive linear combinations, i.e. for any finite subset $S \subseteq C$ and $\lambda : S \to \mathbb{R}_{\geq 0}$, we have $\sum_{s \in S} \lambda_s s \in C$. The cone generated by $S \subseteq \mathbb{R}^n$ is the smallest cone that contains $S$ and is defined via

$$\text{cone}(S) := \left\{ \sum_{i=1}^{k} \lambda_i x_i \mid k \in \mathbb{N}, x_i \in S, \lambda_i \geq 0 \text{ for } i = 1, ..., k \right\}.$$

We formally define $\text{cone}(\emptyset) = \{\mathbf{0}\}$. Clearly $\text{cone}(S)$ is a convex cone.

**Remark 2** (Addition of cones): Let $A, B \subseteq \mathbb{R}^n$. Then $\text{cone}(A) + \text{cone}(B) = \text{cone}(A \cup B)$.

**Remark 3**: If a cone $C$ has a finite generator, then it also has a minimum generator.

**Definition 4** (Generator, Conically Independent Generator, Minimum Generator): A set $S' \subseteq S$ is called generator of $\text{cone}(S)$ if $\text{cone}(S') = \text{cone}(S)$. When $\text{cone}(S)$ is clear from context, we may also just say that $S'$ is a generator.

$S'$ is called conically independent generator if $S'$ is a generator and there is no $S'' \subsetneq S'$ that is a generator. $S'$ is called minimum generator if for any generator $S''$ it holds that $|S'| \leq |S''|$.

The following optimality invariant follows directly.

**Lemma 5**: Let $S$ be a conically independent (minimum) generator of $\text{cone}(S)$. Let $S' \subseteq S$. Then $S'$ is a conically independent (minimum) generator of $\text{cone}(S')$.

*Proof*: Conically independent: Suppose $S'$ is not conically independent. Then there is a generator $S'' \subsetneq S'$ of $\text{cone}(S')$. Now $\text{cone}(S) = \text{cone}(S \setminus S' \cup S'')$ due to Remark 2, a contradiction.

Minimum: Suppose $S'$ is not minimum. Then there is a minimum generator $S''$ of $\text{cone}(S')$ with $|S''| < |S'|$. Now $\text{cone}(S) = \text{cone}(S \setminus S' \cup S'')$ due to Remark 2, a contradiction. □

## 3 Conically independent and minimum generators

The main idea of the proof is to recognize the degenerated case when $\text{cone}(S)$ contains a non-trivial linear subspace, as the situation is nice when $\text{cone}(S)$ is pointed:

**Definition 6** (Lineality space, Pointed): The greatest linear subspace $L$ contained in $\text{cone}(S)$ is called lineality space of $\text{cone}(S)$. $\text{cone}(S)$ is called pointed if $L = \{0\}$.

**Remark 7**: It holds that $L = \text{cone}(S) \cap \text{cone}(-S) = \text{cone}(S) \cap -\text{cone}(S)$.

**Theorem 8**: If $\text{cone}(S)$ is pointed, then any conically independent generator is a minimum generator.

*Proof*: Follows directly from 2.10.2 in [3]. □

**Remark 9**: If $S \subseteq \mathbb{R}^n_{\geq 0}$, then cone($S$) is pointed: It holds that cone($S$) $\subseteq \mathbb{R}^n_{\geq 0}$ such that cone($-S$) $\subseteq \mathbb{R}^{n(\leq 0)}$ implying cone($S$) $\cap$ cone($-S$) = $\{0\}$.

When $L \neq \{0\}$, we split up cone($S$) into its "degenerated" lineal part and into its conical part:

**Definition 10** (Lineal part, Conical part): Let $S'$ be a generator of cone($S$). We define its lineal part $S'_L = S' \cap L$, and its conical part $S'_C = S' \setminus S'_L$.

We collect some properties about this decomposition. Let $\pi_{L^\perp} : \mathbb{R}^n \to L^\perp$ denote the orthogonal projection on the subspace $L^\perp$.

**Lemma 11**:
1) $L = \text{cone}(S_L)$
2) $\text{cone}(S) = L + \text{cone}(S_C) = L + \text{cone}(\pi_{L^\perp}(S_C))$

*Proof*:
1) See [2].
2) For the first equality note that $S = S_L \cup S_C$. Use 1) and Remark 2. For the second equality, set $l_s = s - \pi_{L^\perp}(s) \in L$ for $s \in S_C$. Take note of the following calculation.
$$l + \sum_{s \in S_C} \lambda_s s = l + \sum_{s \in S_C} \lambda_s (s - \pi_{L^\perp}(s) + \pi_{L^\perp}(s))$$
$$= l + \sum_{s \in S_C} \lambda_s l_s + \sum_{s \in S_c} \lambda_s \pi_{L^\perp}(s)$$

This already shows "$\subseteq$". For "$\supseteq$", subtract $\sum_{s \in S_C} \lambda_s l_s$ from both sides and note that $-\sum_{s \in S_C} \lambda_s l_s \in \text{cone}(L)$. □

By Lemma 11.2) we can use $\text{cone}(\pi_{L^\perp}(S_C))$ instead of $\text{cone}(S_C)$. The nice thing about this is that when projecting onto $L^\perp$, it turns out that cone($S_C$) stays a cone and is guaranteed to be pointed. Further, $\text{cone}(\pi_{L^\perp}(S'_C))$ is the same for any generator $S'$ of cone($S$).

**Lemma 12**:
1) $\pi_{L^\perp}(\text{cone}(S)) = \pi_{L^\perp}(\text{cone}(S_C)) = \text{cone}(\pi_{L^\perp}(S_C))$
2) If $S'$, $S''$ are generators of cone($S$), then $\text{cone}(\pi_{L^\perp}(S'_C)) = \text{cone}(\pi_{L^\perp}(S''_C))$.
3) $\text{cone}(\pi_{L^\perp}(S_C))$ is pointed.

*Proof*:
1) Both equalities follow directly by the linearity of projections:
   The first equality is shown by $\pi_{L^\perp}(\text{cone}(S)) = \pi_{L^\perp}(L + \text{cone}(S_C)) = \pi_{L^\perp}(\text{cone}(S_C))$. The second equality follows from
   $$x \in \pi_{L^\perp}(\text{cone}(S_C)) \iff \exists \lambda \in \mathbb{R}^{S_C}_{\geq 0} : x = \pi_{L^\perp}\left(\sum_{s \in S_C} \lambda_s s\right)$$
   $$\iff \exists \lambda \in \mathbb{R}^{S_C}_{\geq 0} : x = \sum_{s \in S_C} \lambda_s \pi_{L^\perp}(s)$$
   $$\iff x \in \text{cone}(\pi_{L^\perp}(S_C)).$$
2) By 1) it follows that $\text{cone}(\pi_{L^\perp}(S'_C)) = \pi_{L^\perp}(\text{cone}(S')) = \pi_{L^\perp}(\text{cone}(S)) = \pi_{L^\perp}(\text{cone}(S'')) = \text{cone}(\pi_{L^\perp}(S''_C))$
3) Let $s \in \pi_{L^\perp}(\text{cone}(S)) \cap -\pi_{L^\perp}(\text{cone}(S))$. This implies $s + l_1 \in \text{cone}(S)$ and $s + l_2 \in \text{cone}(-S)$ for some $s \in L^\perp, l_1, l_2 \in L$. Then also $s \in \text{cone}(S) \cap \text{cone}(-S) = L$. Since $s \in L^\perp \cap L$ it follows that $s = 0$. □

Since $\pi_{L^\perp}(\text{cone}(S))$ is pointed, we can prove one of the two key ingredients for proving Theorem 16.

**Lemma 13**: Suppose $S$ is a conically independent generator. Then $\pi_{L^\perp}(S_C)$ is a minimum generator of $\pi_{L^\perp}(\text{cone}(S))$ and $|\pi_{L^\perp}(S_C)| = |S_C|$.

*Proof*: We first prove that $|\pi_{L^\perp}(S_C)| = |S_C|$. It suffices to show $|\pi_{L^\perp}(S_C)| \geq |S_C|$. Suppose there are $s, s' \in S_C$ s.t. $\pi_{L^\perp}(s) = \pi_{L^\perp}(s')$. Then $s - s' \in \pi_{L^\perp}^{-1}(\{0\}) = L = \text{cone}(S_L)$ and therefore $s = s' + l$ for some $l \in L$, implying $\text{cone}(S \setminus \{s\}) = \text{cone}(S)$, a contradiction.

By Lemma 12.1) we know that $\pi_{L^\perp}(\text{cone}(S)) = \text{cone}(\pi_{L^\perp}(S_C))$. Suppose that there is $s \in \pi_{L^\perp}(S_C)$ s.t. $\text{cone}(\pi_{L^\perp}(S_C)) = \text{cone}(\pi_{L^\perp}(S_C) \setminus \{s\})$. Let $s_c \in S_C$ be the uniquely defined element such that $\pi_{L^\perp}(s_c) = s$. Then $\text{cone}(S) \supseteq \text{cone}(S \setminus \{s_c\}) \supseteq \text{cone}(S_L) + \text{cone}(S_C \setminus \{s_c\}) = L + \text{cone}(\pi_{L^\perp}(S_C \setminus \{s_c\})) = L + \text{cone}(\pi_{L^\perp}(S_C)) = L + \text{cone}(S_C) = \text{cone}(S)$. This shows that $\pi_{L^\perp}(S_C)$ is a conically independent generator. By Lemma 12.3) we know that $\text{cone}(\pi_{L^\perp}(S_C))$ is pointed. Use Theorem 8 on $\text{cone}(\pi_{L^\perp}(S_C))$. □

Before we can put things together, we need to prove the statement also for the lineal part. It turns out that this is where our approximation factor of 2 appears.

**Lemma 14**: Let $d := \dim(L)$ be the dimension of L and $V := \{v_1, ..., v_d\} \subseteq L$ be a set of $d$ linearly independent vectors. Define $V^- := \left\{\sum_{i=1}^d \lambda_i v_i \mid \lambda_1, ..., \lambda_d < 0\right\}$. For any $A \subseteq L$ we have that
1) $\text{cone}(V \cup A) = L \iff \text{cone}(A) \cap V^- \neq \emptyset$.
2) If $V \cup A$ is a conically independent generator of $L$, then $A =: \{a_1, ..., a_k\}$ is finite and $\forall a \in \text{cone}(A) \cap V^- : \left(\exists \lambda_1, ..., \lambda_k \geq 0 : a = \sum_{i=1}^k \lambda_i a_i \implies \lambda_1, ..., \lambda_k > 0\right)$.
3) If $V \cup A$ is a conically independent generator of $L$, then $(a_1, ..., a_k)$ is linearly independent.

*Proof*:
1) "$\impliedby$": Let $w \in \text{cone}(A) \cap V^-$. Then $w = \sum_{i=1}^d \lambda_i v_i$ for some $\lambda_1, ..., \lambda_d < 0$. Let $l = \sum_{i=1}^d l_i v_i \in L$, $l_i \in \mathbb{R}$. Let $n \in \mathbb{N}$ such that $n\lambda_i < l_i$ for all $i = 1, ..., d$. Then $l = nw + \sum_{i=1}^d (l_i - n\lambda_i) v_i$.
"$\implies$": We know that $x = -\sum_{i=1}^d v_i \in V^- \cap L$. In particular, since $x \in L$, $x = v + a$ for some $v \in \text{cone}(V), a \in \text{cone}(A)$. It follows that $a = x - v \in V^-$.
2) Let $a \in \text{cone}(A) \cap V^-$ with $a = \sum_{i=1}^k \lambda_i a_i$ for some $\lambda_1, ..., \lambda_k \geq 0$ and pairwise different $a_1, ..., a_k \in A$. If $A \neq \{a_1, ..., a_k\}$, then there is $z \in A \setminus \{a_1, ..., a_k\}$. Then $a \in \text{cone}(A \setminus \{z\})$, which by 1) proves that $\text{cone}(V \cup A \setminus \{z\}) = L$, which is a contradiction to the conic independence of $V \cup A$. The same argument also shows that $\lambda_i > 0$ for all $i \in [k]$.
3) Let $\boldsymbol{A} \in \mathbb{R}^{n \times k} := (a_1 \ a_2 \ ... \ a_k)$ and assume $(a_1, ..., a_k)$ are linearly dependent. This implies the existence of $z \neq 0 \in \mathbb{R}^k$ such that $\boldsymbol{A}z = 0$. By 1) there exists $x^* \geq 0$ s.t. $\boldsymbol{A}x^* \in \text{cone}(A) \cap V^-$. W.l.o.g. $z$ has at least one negative component (If $z \geq 0$, replace $z$ by $-z$). Set $m := \max\left\{l \in \mathbb{R}_{\geq 0} \mid (x^* + lz)_i \geq 0 \text{ for all } i \in [k]\right\}$ and note that the maximum is well defined, as the set is non-empty, closed and bounded due to $z$ having a negative component. Then $(x^* + mz) \geq 0$ and by maximality of $m$ also $(x^* + mz)_i = 0$ for some $i \in [k]$. Further $\boldsymbol{A}(x^* + mz) = \boldsymbol{A}x^* \in \text{cone}(A) \cap V^-$, in contradiction to 2). □

**Theorem 15**: Let $G \subseteq \mathbb{R}^n$ be a minimum generator of $\text{cone}(S)$ and $H \subseteq \mathbb{R}^n$ be any conically independent generator of $L$. Then $|H| \leq 2|G|$

*Proof*: Let $d := \dim(L)$ be the dimension of L. Clearly, $L = \text{cone}(G) \subseteq \text{span}(G) \subseteq L$. Therefore $G$ contains at least $d$ linearly independent vectors. It follows that $|G| \geq d$, and analogously $H$ contains at least $d$ linearly independent vectors and $|H| \geq d$. We show $|H| \leq 2d$. Denote by $v_1, ..., v_d \in H$ linearly independent vectors and let $H =: V \cup A := \{v_1, ..., v_d\} \cup A$ where $A = H \setminus V$. Due to Lemma 14.2) and Lemma 14.3), we get that $A$ is finite and linearly independent, and can conclude $|A| \leq d$. □

**Theorem 16**: If $S'$ is a conically independent generator of $\text{cone}(S)$ and $S''$ is a minimum generator of $\text{cone}(S)$, then $|S'| \leq 2\,|S''|$

*Proof*: By Lemma 5 and Lemma 11.1), $S'_L$ is a conically independent generator of $L$ and $S''_L$ is minimum generator of $L$. By Theorem 15 it follows that $|S'_L| \leq 2|S''_L|$. We know due to Lemma 13 that both $\pi_{L^\perp}(S'_C)$ and $\pi_{L^\perp}(S''_C)$ are a minimum generator of their respective conic hull, and by Lemma 12.2) $\text{cone}(\pi_{L^\perp}(S'_C)) = \text{cone}(\pi_{L^\perp}(S''_C))$. This shows that $|\pi_{L^\perp}(S'_C)| = |\pi_{L^\perp}(S''_C)|$, and by Lemma 13 we get $|S'_C| = |S''_C|$. In total we get:

$|S'| = |S'_L| + |S'_C| \leq 2|S''_L| + |S'_C| \leq 2(|S''_L| + |S''_C|) = 2\,|S''|$ □

## 4 Algorithms

**Lemma 17**: Linear programs can be solved in polynomial time $O(\tilde{n}^{3.5}L)$, where $L$ is the number of bits needed to input the linear program and $\tilde{n}$ is the number of variables. Therefore linear programs can be solved in polynomial-time (in number of variables/input size).

*Proof*: This is Karmarkar's algorithm, see [5]. □

**Lemma 18**: Let $c \in \mathbb{R}^n$. The property $c \in \text{cone}(S)$ can be checked in polynomial time in $n$ and $|S|$.

*Proof*: Let $S =: \{s^{(1)}, ..., s^{(|S|)}\}$. Then $c \in \text{cone}(S)$ if and only if

$$\sum_{i=1}^{k} y_i s^{(i)} = c$$
$$y \in \mathbb{R}_{\geq 0}^{|S|}$$
(1)

This is an LP without objective function. We can add the constant 0 function as objective, and use Lemma 17. □

**Algorithm 19**: We state a polynomial-time algorithm to find a conically independent generator given a finite generator $S =: \{s^{(1)}, ..., s^{(|S|)}\}$.

1: **function** FIND-CONICALLY-INDEPENDENT-GENERATOR(S)
2:     $S' \leftarrow S$
3:     **for** $i \in \{1, ..., |S|\}$ **do**
4:         **if** $s^{(i)} \in \text{cone}(S' \setminus \{s^{(i)}\})$ **then**
5:             $S' \leftarrow S' \setminus \{s^{(i)}\}$
6:     **return** $S'$

*Proof*: Correctness: Clearly, $S'$ remains a generator of $\text{cone}(S)$ in each step, since we explicitly check this condition in the if statement. Furthermore, $S'$ is conically independent: Suppose $s^{(i)} \in \text{cone}(S' \setminus \{s^{(i)}\})$ for some $s^{(i)} \in S'$. Denote by $S'_i$ the value of $S'$ during the loop in iteration $i$. Clearly, $S'_i \supseteq S'$, which implies $s^{(i)} \in \text{cone}(S'_i \setminus \{s^{(i)}\})$. But then Line 5 is executed in iteration $i$, a contradiction.

Runtime: Line 4 can be executed in polynomial time, see Lemma 18. It is executed $|S|$ times, thus the runtime stays polynomial. □

Note that in [2], there is a specialized simplex-based algorithm for finding conically independent generators.

**Algorithm 20**: Let $S := \{s^{(1)}, ..., s^{(|S|)}\}$. We state our polynomial time procedure to find a minimum cardinality generator of cone($S$).

1:   **function** FIND-MINIMUM-CARDINALITY-GENERATOR(S)
2:       $S' \leftarrow$ FIND-CONICALLY-INDEPENDENT-GENERATOR($S$)
3:       $S'_L \leftarrow \{s \in S' \mid -s \in \text{cone}(S')\}$
4:       $B \leftarrow$ basis of $L$ in $S'_L$
5:       **return** $B \cup \left\{-\sum_{x \in B} x\right\} \cup S' \setminus S'_L$

*Proof*: Correctness: $S'$ in Line 2 is conically independent due to Algorithm 19. $S'_L$ is determined correctly, see Definition 10. $B \cup \left\{-\sum_{x \in B} x\right\}$ is a minimum generator of $L$, see Lemma 14. The rest follows as in Theorem 16.

Runtime: Line 2 is executed in polynomial time, see Algorithm 19. Line 3 is executed in polynomial time, as it is an $|S'|$-fold application of Lemma 18. Line 4 can be realized in polynomial time, since we can use the Gauss algorithm to find a basis. □